\documentclass[11pt,reqno,a4paper]{amsart}
\usepackage{euscript}
\usepackage{comment}

\theoremstyle{plain}
\numberwithin{equation}{section}
\newtheorem{theorem}{Theorem}
\newtheorem{proposition}{Proposition}
\newtheorem{lemma}{Lemma}
\newtheorem{definition}{Definition}
\newtheorem{corollary}{Corollary}
\theoremstyle{remark}
\newtheorem{remark}{Remark}
\newtheorem{example}{Example}
\renewcommand{\epsilon}{\varepsilon}
\renewcommand{\phi}{\varphi}
\DeclareMathOperator{\Ima}{Im}
\DeclareMathOperator{\Ker}{Ker}

\DeclareMathOperator{\codim}{codim}

\def\Id{\text{\rm Id}}

\usepackage{xcolor}
\usepackage{amsthm}
\usepackage{amsmath}
\usepackage{amsfonts}

\begin{document}
\title[$h$-dichotomies for noninvertible evolution families]{$h$-dichotomies via  uniform noncriticality and expansiveness for  noninvertible evolution families}
\author[Davor Dragi\v cevi\' c]{Davor Dragi\v cevi\' c } 
\address{Department of Mathematics, University of Rijeka, Croatia}
\email[Davor Dragi\v cevi\' c]{ddragicevic@math.uniri.hr}

\begin{abstract}
We establish characterizations of  $h$-dichotomies for noninvertible evolution families acting on Banach spaces, in terms of appropriate notions of expansiveness and uniform noncriticality. Our results extend the corresponding characterizations we previously obtained for invertible evolution families ((Math. Ann. 395 (2026), 17 pp.). To the best of our knowledge, no such characterizations have previously been established for general noninvertible evolution families. A key difference is that noninvertibility makes the structure of the stable spaces more delicate: in our positive-half-line result, finite codimension must be imposed at every time, rather than only at the initial time. We show by an example that the latter, weaker condition is insufficient. We also study evolution families on the negative half-line, using duality to obtain corresponding characterizations in terms of dual expansiveness and dual uniform noncriticality. General growth rates are treated by a time-rescaling argument.
\end{abstract}
\maketitle
\section{Introduction}
The notion of an exponential dichotomy introduced by Perron~\cite{Per} plays a central role in the qualitative theory of nonautonomous differential equations and dynamical systems. It expresses the existence of an invariant splitting into stable and unstable directions, with exponential contraction in forward time on the stable component and exponential contraction in backward time on the unstable component. The classical theory, originating in the work of Perron, has been developed extensively; see, for example, \cite{CL, Coppel,DK, Henry, Pot, SY-book-2002}. For evolution families that are not invertible, the definition of a dichotomy requires particular care: backward evolution is available only on the unstable subspaces, on which the forward operators are required to induce isomorphisms.

Exponential rates are not always appropriate for describing nonautonomous dynamics. This motivates the study of generalized dichotomies whose contraction and expansion rates are prescribed by a growth function $h$. Such $h$-dichotomies encompass exponential dichotomies and allow for a substantially wider range of asymptotic behavior. For the origins of this approach we refer to~\cite{Martin, M, NP}.

An important problem is to characterize dichotomies without assuming an invariant splitting in advance. In the classical exponential setting, this problem was investigated through notions of noncriticality and expansiveness, notably by Coppel~\cite{Coppel} and Palmer~\cite{Palmer}. More recently, Elorreaga, Pe\~na and Robledo~\cite{EPR, EPR2} extended this perspective to $h$-dichotomies for finite-dimensional nonautonomous differential equations (see~\cite{WXZ} for related earlier work). In \cite{D}, the corresponding results of~\cite{EPR} were extended to invertible evolution families on Banach spaces by means of a time-rescaling approach introduced in a systematic manner by Dragi\v cevi\' c and Silva~\cite{DS, DS1} (building on the previous work~\cite{DSS} on polynomial dichotomies). In particular, under suitable growth assumptions and a finite-dimensional complement condition, $h$-dichotomy was characterized by $h$-expansiveness and uniform $h$-noncriticality.

The purpose of the present paper is to develop analogous characterizations for \emph{noninvertible} evolution families on Banach spaces. To the best of our knowledge, the existing characterizations in terms of expansiveness and noncriticality do not cover this general noninvertible setting. This extension is not merely formal. The definitions of expansiveness and noncriticality must be expressed using forward evolution, since the operators $T(t,s)$ need not admit inverses on the whole space. More importantly, noninvertibility changes the geometry of the spaces of initial conditions with bounded forward orbits, and this affects the hypotheses needed for a characterization.

Our first main result concerns evolution families on the positive half-line. For each $s\geq 0$, consider the subspace
\[
\mathcal S(s)=\left\{v\in X:\sup_{t\geq s}\|T(t,s)v\|<\infty\right\}.
\]
We prove that, under bounded growth and the assumption that $\operatorname{codim}\mathcal S(s)$ is the same finite number for every $s\geq0$, exponential dichotomy, exponential expansiveness and uniform noncriticality are equivalent (Theorem~\ref{T1}). Time rescaling then yields the analogous characterization of $h$-dichotomies (Theorem~\ref{T2}).

A significant feature of the noninvertible case is that the finite-codimension hypothesis is more demanding than in the invertible setting. In \cite[Theorem~1]{D}, it suffices to assume that $\mathcal S(0)$ has finite codimension: invertibility propagates this property, with the same codimension, to every later time. For noninvertible evolution families, this propagation is no longer automatic. As explained in Remark~6 and demonstrated by Example~1, finite codimension of $\mathcal S(0)$ alone does not suffice for our characterization. Thus the additional hypothesis in Theorem~\ref{T1} reflects a genuine obstruction, rather than a technical artifact of its proof.

We also investigate evolution families on the negative half-line, where a direct analogue of the positive-half-line characterization is not available in the same form. By considering the time-reversed adjoint evolution family (on a Hilbert space), we obtain a duality principle for dichotomies and formulate corresponding criteria in terms of dual $h$-expansiveness and dual uniform $h$-noncriticality. This approach highlights the asymmetry between forward and backward dynamics when the evolution operators are not invertible.

The paper is organized as follows. We first recall the relevant definitions and establish the time-rescaling framework. We then prove the characterizations on the positive half-line, discuss the necessity of the strengthened codimension assumption, and provide a counterexample. Finally, we study the negative half-line and derive the dual characterizations.

\section{Preliminaries}
Throughout this paper $X=(X, \|\cdot \|)$ will be a Banach space. By $\mathcal B(X)$ we denote the space of all bounded linear operators on $X$ equipped with the operator norm, which we also denote by $\| \cdot \|$.

We begin by introducing the notion of a growth rate.
\begin{definition}
Let $J\in \{\mathbb R, \mathbb R^+, \mathbb R^-\}$. A growth rate $h \colon J \to (0, \infty)$ is  any increasing and continuous map which is a bijection onto its image and satisfies the following:
\begin{itemize}
\item $h(0)=1$;
\item $\lim\limits_{t\to \infty}h(t)=\infty$ if $J\in \{\mathbb R, \mathbb R^+\}$ and $\lim\limits_{t\to -\infty}h(t)=0$ if $J\in \{\mathbb R, \mathbb R^-\}$.
\end{itemize}
\end{definition}

Next, we recall the notion of an evolution family.
\begin{definition}\label{EF}
    Let $J\in \{\mathbb R, \mathbb R^+, \mathbb R^-\}$. We say that $\mathcal T=\{T(t, s): \ t\ge s, \ t, s\in J\}$ is an \emph{evolution family} over $J$ if the following holds:
    \begin{itemize}
\item $T(t, t)=\Id$ for $t\in J$, where $\Id$ denotes the identity operator on $X$;
\item for $t, s, r\in J$ with $t\ge s\ge r$,
\[
T(t, s)T(s, r)=T(t, r);
\]
\item for $s\in J$ and $v\in X$, the map $t\mapsto T(t, s)v$ is continuous on $[s, \infty)\cap J$, and for each $t\in J$ and $v\in X$, the map $s\mapsto T(t, s)v$ is continuous on $J\cap (-\infty, t]$.
    \end{itemize}
    If in addition,  $T(t, s)$ is an invertible operator for each $t\ge s$, we say that $\mathcal T$ is an \emph{invertible} evolution family.
\end{definition}

\begin{remark}
If  $\mathcal T=\{T(t, s): \ t\ge s, \ t, s\in J\}$ is an invertible evolution family over $J\in \{\mathbb R, \mathbb R^+, \mathbb R^-\}$, we can define
\[
T(t, s):=T(s, t)^{-1}, \quad \text{for $t, s\in J$, $t<s$.}
\]
In this case, we can write $\mathcal T=\{T(t, s): \ t, s\in J\}$.

\end{remark}

\begin{definition}
Let  $J\in \{\mathbb R, \mathbb R^+, \mathbb R^-\}$, $h\colon J\to (0, \infty)$ be a growth rate and $\mathcal T=\{T(t, s): \ t\ge s, \ t, s\in J\}\subset \mathcal B(X)$ an evolution family. We say that $\mathcal T$ has \emph{$h$-bounded growth} if there are $K, \omega>0$ such that 
\begin{equation}\label{HBG}
\|T(t, s)\| \le K \left (\frac{h(t)}{h(s)}\right)^\omega, \quad \text{for $t, s\in J$ with $t\ge s$.}
\end{equation}

\end{definition}

\begin{remark}
In the particular case when $h(t)=e^t$, we say that $\mathcal T$ has bounded growth.
\end{remark}

We recall the notion of $h$-dichotomy.
\begin{definition}
Let  $J\in \{\mathbb R, \mathbb R^+, \mathbb R^-\}$, $h\colon J\to (0, \infty)$ be a growth rate and $\mathcal T=\{T(t, s): \ t\ge s, \ t, s\in J\}\subset \mathcal B(X)$ an evolution family. We say that $\mathcal T$ admits an \emph{$h$-dichotomy}  if there is a family $\{P(t): \ t\in J\}\subset \mathcal B(X)$ of projections on $X$ and constants $D, \lambda >0$ such that:
\begin{enumerate}
\item for $t, s\in J$ with $t\ge s$,
\begin{equation}\label{inv0}
P(t)T(t, s)=T(t, s)P(s),
\end{equation}
and \begin{equation}\label{inv}T(t, s)\rvert_{\Ker P(s)}\colon \Ker P(s)\to \Ker P(t) \quad  \text{is invertible};\end{equation}

\item for $t, s\in J$ with $t\ge s$,
\begin{equation}\label{HD1}
\|T(t, s)P(s)\| \le D \left (\frac{h(t)}{h(s)}\right )^{-\lambda};
\end{equation}
\item for $t, s\in J$ with $t\le s$,
\begin{equation}\label{HD2}
\|T(t, s)(\Id-P(s))\| \le D\left (\frac{h(s)}{h(t)}\right)^{-\lambda},
\end{equation}
where
\[
T(t, s):=\left (T(s, t)\rvert_{\Ker P(t)}\right )^{-1}\colon \Ker P(s)\to \Ker P(t).
\]
\end{enumerate}
\end{definition}
\begin{remark}
In the particular case when $h(t)=e^t$, we recover the notion of an \emph{exponential dichotomy}.
\end{remark}
We introduce the notion of $h$-expansivity for not necessarily invertible evolution families.
\begin{definition}
Let  $J\in \{\mathbb R, \mathbb R^+, \mathbb R^-\}$, $h\colon J\to (0, \infty)$ be a growth rate and $\mathcal T=\{T(t, s): \ t\ge s, \ t, s\in J\}\subset \mathcal B(X)$ an evolution family. We say that $\mathcal T$ is \emph{$h$-expansive} if there exist $L, \beta>0$ such that 
\begin{equation}\label{HEXP}
\|T(t, a)v\| \le L\left ( \left (\frac{h(t)}{h(a)}\right )^{-\beta}\|v\|+ \left (\frac{h(b)}{h(t)}\right )^{-\beta}\|T(b, a)v\| \right ),
\end{equation}
for $v\in X$ and $t\in [a, b] \subset J$.
\end{definition}

\begin{remark}
In the particular case when $h(t)=e^t$, we say that $\mathcal T$ is exponentially expansive.
\end{remark}

The following simple result shows that for invertible evolution families over $J=\mathbb R^+$, the concept of $h$-expansiveness coincides with the one introduced in~\cite[Definition 6]{D} (adapting the definition of Palmer~\cite{Palmer} for ODEs).
\begin{proposition}
Let  $J\in \{\mathbb R, \mathbb R^+, \mathbb R^-\}$,  $h\colon J\to (0, \infty)$ be a growth rate and $\mathcal T=\{T(t, s):  \ t, s\in J\}\subset \mathcal B(X)$ an invertible evolution family. The following properties are equivalent:
\begin{enumerate}
\item[(a)] $\mathcal T$ is $h$-expansive;
\item[(b)] there exist $L, \beta>0$ such that
\begin{equation}\label{HEXPP}
\|v\| \le L\left ( \left (\frac{h(t)}{h(a)}\right )^{-\beta}\|T(a, t)v\|+\left (\frac{h(b)}{h(t)}\right )^{-\beta}\|T(b, t)v\| \right ),
\end{equation}
for $v\in X$ and $t\in [a, b]\subset J$.
\end{enumerate}

\end{proposition}

\begin{proof}
$(a)\implies (b)$ Suppose that $\mathcal T$ is $h$-expansive and let $L, \beta>0$ be such that~\eqref{HEXP} holds for $v\in X$ and $t\in [a, b]\subset J$. 

Take an arbitrary $v\in X$ and $t\in [a, b]\subset J$. By applying~\eqref{HEXP} to $T(a, t)v$ instead of $v$, we conclude that~\eqref{HEXPP} holds.

$(b) \implies (a)$ This implication follows by analogous reasoning. Namely, applying~\eqref{HEXPP} to $T(t, a)v$ instead of $v$, we deduce that~\eqref{HEXP} holds.

\end{proof}

Finally, we introduce the notion of uniform $h$-noncriticality for not necessarily invertible evolution families.
\begin{definition}
Let  $J\in \{\mathbb R, \mathbb R^+, \mathbb R^-\}$, $h\colon J\to (0, \infty)$ be a growth rate and $\mathcal T=\{T(t, s): \ t\ge s, \ t, s\in J\}\subset \mathcal B(X)$ an evolution family. We say that $\mathcal T$ is uniformly \emph{$h$-noncritical} if there exist $\theta \in (0, 1)$ and $C>0$ such that 
\begin{equation}\label{HNONC}
\|T(h^{-1}(e^C h(t)), t)v\| \le \theta \sup \{\|T(u, t)v\|: \ u\ge t, \ \ln h(u)-\ln h(t)\le 2C\},
\end{equation}
for all $v\in X$ and $t\in J$ such that $[t, h^{-1}(e^{2C}h(t))]\subset J $.

\end{definition}

\begin{remark}
In the particular case when $h(t)=e^t$, we say that $\mathcal T$ is uniformly noncritical.
\end{remark}
The following proposition shows that for invertible evolution families over $J=\mathbb R^+$, the concept of uniform $h$-noncriticality coincides with the one introduced in~\cite[Definition 7]{D} (adapting the notion introduced for ODE's in~\cite{K, MS}).

\begin{proposition}
Let  $J\in \{\mathbb R, \mathbb R^+, \mathbb R^-\}$, $h\colon J\to (0, \infty)$ be a growth rate and $\mathcal T=\{T(t, s): \ t, s\in J\}\subset \mathcal B(X)$ an invertible evolution family. The following properties are equivalent:
\begin{enumerate}
\item[(a)] $\mathcal T$ is uniformly $h$-noncritical;
\item[(b)] there exist $\theta \in (0, 1)$ and $C>0$ such that 
\begin{equation}\label{H616}
\| v\| \le \theta \sup \{\|T(u, t)v\|: \ |\ln h(u)-\ln h(t)| \le C\},
\end{equation}
for all $v\in X$ and $t\in J$ such that $[h^{-1}(e^{-C} h(t)), h^{-1}(e^C h(t))]\subset J$.
\end{enumerate}

\end{proposition}

\begin{proof}
$(a)\implies (b)$ Suppose that $\mathcal T$ is uniformly $h$-noncritical and let $\theta \in (0, 1)$ and $C>0$ be such that~\eqref{HNONC} holds for $v\in X$ and $t\in J$ such that $[t, h^{-1}(e^{2C}h(t))]\subset J $.

Take now $v\in X$ and $t\in J$ such that $[h^{-1}(e^{-C} h(t)), h^{-1}(e^C h(t))]\subset J$. By applying~\eqref{HNONC} for $h^{-1}(e^{-C}h(t))$ instead of $t$ (which we can since $[h^{-1}(e^{-C}h(t)), h^{-1}(e^C h(t))]\subset J$), we have 
\[
\begin{split}
&\|T(t, h^{-1}(e^{-C}h(t)))w\|\\
&\le  \theta \sup \{\|T(u, h^{-1}(e^{-C}h(t)))w\|: \ |\ln h(u)-\ln h(t)| \le C\}.
\end{split}
\]
By substituting $T(h^{-1}(e^{-C} h(t)), t)v$ instead of $w$, we obtain~\eqref{H616}. Thus, $(b)$ holds.

$(b)\implies (a)$ Take $v\in X$ and $t\in J$ such that $[t, h^{-1}(e^{2C}h(t))]\subset J $. By applying~\eqref{H616} for $h^{-1}(e^C h(t))$ instead of $t$, we have 
\[
\|w\| \le \theta \sup \{\|T(u, h^{-1}(e^{C}h(t)))w\|: \ u\ge  t, \ \ln h(u)-\ln h(t)\le 2C\},
\]
for $w\in X$. Taking $w=T(h^{-1}(e^C h(t)), t)v$ yields~\eqref{HNONC}.

\end{proof}

\section{Time-rescaling}
We recall the construction introduced in~\cite{D} (see also~\cite{DS1}). The proof of the following result was established in~\cite[Proposition 2]{D} in the case when $J=\mathbb R^+$. The proof in the case when $J\in \{\mathbb R, \mathbb R^-\}$ is completely analogous.

\begin{theorem}\label{TR}
Let  $J\in \{\mathbb R, \mathbb R^+, \mathbb R^-\}$, $h\colon J\to (0, \infty)$ be a growth rate and $\mathcal T=\{T(t, s): \ t\ge s, \ t, s\in J\}\subset \mathcal B(X)$ an evolution family. Set
\[
T_h(t, s):=T(h^{-1}(e^t), h^{-1}(e^s)), \quad t, s\in J, \ t\ge s.
\]
Then $\mathcal T^h=\{T_h(t, s): \ t, s \in J, \ t\ge s\}$ is an evolution family. Moreover, $\mathcal T$ admits an $h$-dichotomy if and only if $\mathcal T^h$ admits an exponential dichotomy.\end{theorem}

The next two results are similar in nature to Theorem~\ref{TR}, showing that $h$-expansivity (resp. uniform $h$-noncriticality) of $\mathcal T$ is equivalent to exponential expansivity (resp. uniform noncriticality) of $\mathcal T^h$.
\begin{proposition}\label{TR1}
    Let  $J\in \{\mathbb R, \mathbb R^+, \mathbb R^-\}$, $h\colon J\to (0, \infty)$ be a growth rate and $\mathcal T=\{T(t, s): \ t\ge s, \ t, s\in J\}\subset \mathcal B(X)$ an evolution family. Then $\mathcal T$ is $h$-expansive if and only if $\mathcal T^h$ is exponentially expansive.

\end{proposition}

\begin{proof}
Assume that $\mathcal T$ is $h$-expansive and let $L, \beta>0$ be such that~\eqref{HEXP} holds for $t\in [a, b]\subset J$. Take now arbitrary $[a, b]\subset J$ and $t\in [a, b]$. Then, $h^{-1}(e^t)\in [h^{-1}(e^a), h^{-1}(e^b)]\subset J$  and thus using~\eqref{HEXP} we have 
\[
\begin{split}
&\|T_h(t, a)v\| \\
&=\|T(h^{-1}(e^t), h^{-1}(e^a))v\| \\
&\le L\left (\left (\frac{h(h^{-1}(e^t))}{h(h^{-1}(e^a))}\right )^{-\beta}\|v\|+\left (\frac{h(h^{-1}(e^b))}{h(h^{-1}(e^t))}\right)^{-\beta }\|T(h^{-1}(e^b), h^{-1}(e^a))v\|\right ) \\
&= L\left (e^{-\beta (t-a)}\|v\|+e^{-\beta (b-t)}\|T_h(b, a)v\|\right ),
\end{split}
\]
for $v\in X$. Hence, $\mathcal T^h$ is exponentially expansive. The converse implication can be obtained similarly by noting that 
\[
T(t, s)=T_h(\ln h(t), \ln h(s)), \quad t, s\in J, \ t\ge s.
\]

\end{proof}

Similarly to the proof of the previous proposition, one can establish the following.

\begin{proposition}\label{TR2}
    Let  $J\in \{\mathbb R, \mathbb R^+, \mathbb R^-\}$, $h\colon J\to (0, \infty)$ be a growth rate and $\mathcal T=\{T(t, s): \ t\ge s, \ t, s\in J\}\subset \mathcal B(X)$ an evolution family. Then $\mathcal T$ is uniformly $h$-noncritical if and only if $\mathcal T^h$ is uniformly noncritical.

\end{proposition}

\section{Evolution families over $\mathbb R^+$}
For a subspace $Y\subset X$ and $d\in \mathbb N_0$, we write $\codim Y=d$ if there exists a subspace $Z\subset X$ of dimension $d$ such that $X=Y\oplus Z$. Note that $d$ is the dimension of the quotient $X/Y$.

The following is a version of~\cite[Theorem 1]{D} for not necessarily invertible evolution families.  Its proof is inspired by the arguments in~\cite{Coppel, D}.
\begin{theorem}\label{T1}
Let  $J=\mathbb R^+$
and $\mathcal T=\{T(t, s): \ t\ge s, \ t, s\in J\}\subset \mathcal B(X)$ be an evolution family. Furthermore, assume that $\mathcal T$ has  bounded growth and that there exists $d\in \mathbb N_0$ such that \begin{equation}\label{codim}\codim \mathcal S(s)=d, \quad   s\in J,\end{equation} where 
\begin{equation}\label{SS}
\mathcal S(s):=\left \{v\in X: \ \sup_{t\ge s}\|T(t, s)v\|<\infty \right \}.
\end{equation}
Then the following statements are equivalent:
\begin{enumerate}
\item[(a)] $\mathcal T$ admits an exponential dichotomy;
\item[(b)] $\mathcal T$ is exponentially expansive;
\item[(c)] $\mathcal T$ is uniformly noncritical.
\end{enumerate}

\end{theorem}

\begin{proof}
$(a)\implies (b)$ Assume that $\mathcal T$ admits an exponential dichotomy. Take an arbitrary $[a, b]\subset J$. For $t\in [a, b]$ and $v\in X$, we have 
\[
\begin{split}
T(t, a)v &=P(t)T(t, a)v+(\Id-P(t))T(t, a)v \\
&=T(t, a)P(a)v+T(t, b)(\Id-P(b))T(b, a)v,
\end{split}
\]
and consequently using~\eqref{HD1} and~\eqref{HD2} we have 
\[
\|T(t, a)v\| \le D\left (e^{-\lambda (t-a)}\|v\|+e^{-\lambda (b-t)}\|T(b, a)v\|\right ).
\]
Therefore~\eqref{HEXP} holds with $L=D$, $\beta=\lambda$ and $h(t)=e^t$.

$(b)\implies (c)$ Suppose that $\mathcal T$ is exponentially expansive and let $L, \beta>0$ be such that~\eqref{HEXP} holds for  $v\in X$ and $t\in [a, b] \subset J$ with $h(t)=e^t$.

Choose $C>0$ such that $\theta:=2Le^{-\beta C}<1$. Take $v\in X$ and $t\in J$. By applying~\eqref{HEXP} for $a=t$, $b=t+2C$ and $t+C$ instead of $t$, we get
\[
\begin{split}
\|T(t+C, t)v\| &\le L\left (e^{-\beta C}\|v\|+e^{-\beta C}\|T(t+2C, t)v\| \right ) \\
&\le 2Le^{-\beta C}\max \{\|v\|, \|T(t+2C, t)v\|\}\\
&\le \theta \sup \{\|T(u, t)v\|: \ u\in [t, t+2C]\}.
\end{split}
\]
Therefore, \eqref{HNONC} holds (with $h(t)=e^t$).

$(c)\implies (a)$ Since $\mathcal T$ is uniformly noncritical, there exist $\theta \in (0, 1)$ and $C>0$ such that 
\begin{equation}\label{531-1}
\|T(t+C, t)v\| \le \theta \sup \{\|T(u, t)v\|: \ u\in [t, t+2C]\},
\end{equation}
for $t\in J$ and $v\in X$.

Let $Z$ be a subspace of $X$ of dimension $d$ such that 
\begin{equation}\label{SPLIT}
X=\mathcal S(0)\oplus Z.
\end{equation}

 \begin{lemma}
 We have that
\begin{equation}\label{SPLITs}
X=\mathcal S(s)\oplus Z(s), \quad \text{for $s\in J$,}
\end{equation}
where  $Z(s):=T(s, 0)Z$. 

\end{lemma}
\begin{proof}[Proof of the lemma]
Set $Z(s):=T(s, 0)Z$, $s\in J$. We first claim that 
\begin{equation}\label{claim}
T(s, 0)\rvert_Z \colon Z\to Z(s) \quad \text{is invertible for each $s\in J$.}
\end{equation}
Clearly,  $T(s, 0)\rvert_Z$ is onto. Suppose that there is $v\in Z$ such that $T(s, 0)v=0$. Thus, $T(t, 0)v=T(t, s)T(s, 0)v=0$ for $t\ge s$ and consequently $v\in \mathcal S(0)$. We conclude that $v\in \mathcal S(0)\cap Z=\{0\}$, which implies that $T(s, 0)\rvert_Z$ is injective. We have proved~\eqref{claim}. Next, we note that for $s\in J$, $\mathcal S(s)\cap Z(s)=\{0\}$. Indeed, take $w\in \mathcal S(s)\cap Z(s)$, and let $v\in Z$ be such that $w=T(s, 0)v$. Then $v\in \mathcal S(0)\cap Z=\{0\}$, and thus $w=0$. It remains to observe that~\eqref{SPLITs} now follows from $\codim \mathcal S(s)=d$ and $\dim Z(s)=\dim Z=d$.
\end{proof}
For $s<t$, we define $T(s, t)\colon Z(t)\to Z(s)$ to be the inverse of $T(t, s)\rvert_{Z(s)}$.
\begin{lemma}\label{backward}
For $t\in J$ and $v\in Z(t)$, the map $s\mapsto T(s, t)v$ is continuous on $[0, t]$.
\end{lemma}
\begin{proof}[Proof of the lemma]
Choose $w\in Z$ such that $T(t, 0)w=v$. Then $T(s, t)v=T(s, 0)w$, and so the desired continuity follows since $\mathcal T$ is an evolution family (see the last requirement in Definition~\ref{EF}).

\end{proof}

\begin{lemma}\label{LEM:s}
There exist $B, \alpha>0$ such that 
\[
\|T(t, s)v\| \le Be^{-\alpha (t-s)}\|v\|, \quad \text{for $s\in J$, $t\ge s$ and $v\in \mathcal S(s)$.}
\]

\end{lemma}

\begin{proof}[Proof of the lemma]
We first claim  that 
\begin{equation}\label{459new}
\|T(t, s)v\| \le D\|v\|, \quad \text{for $t\ge s\ge 0$ and $v\in \mathcal S(s)$,}
\end{equation}
where $D:=Ke^{\omega C}>0$.
Take $s\in J$ and $v\in \mathcal S(s)\setminus \{0\}$. Then
\[
0<\varrho:=\sup_{t\ge s}\|T(t, s)v\|<+\infty. 
\]
Take $t\ge s+C$. By~\eqref{531-1} (for $t-C$ instead of $t$),
\[
\begin{split}
\|T(t, s)v\| &=\|T(t, t-C)T(t-C, s)v\| \\
&\le \theta \sup  \{\|T(u, t-C)T(t-C, s)v\|: \ u\in [t-C, t+C]\}\\
&=\theta \sup \{\|T(u, s)v\|: \ u\in [t-C, t+C]\}\\
&\le \theta \sup \{\|T(u, s)v\|: \ u\ge s\}\\
&=\theta \varrho <\varrho,
\end{split}
\]
as $\theta \in (0, 1)$.
Hence, 
\[
\varrho=\sup \{\|T(t, s)v\|: \ t\in [s, s+C]\},
\]
and therefore by~\eqref{HBG}, $\varrho \le Ke^{\omega C}\|v\|=D\|v\|$ which yields~\eqref{459new}.

Take now $t\ge s\ge 0$ and $v\in \mathcal S(s)$. Choose $n\in \mathbb N_0$ such that $s+nC \le t<s+(n+1)C$. By~\eqref{531-1}, we have (provided $n\ge 1$) that 
\[
\begin{split}
\|T(t, s)v\| &=\|T(t, t-C) T(t-C, s)v\| \\
&\le \theta \sup \{\|T(u, t-C)T(t-C, s)v\|: \ u \in [t-C, t+C]\}\\
&=\theta \sup \{\|T(u, s)v\|: \ u\in [t-C, t+C]\}.
\end{split}
\]
Furthermore, for $u\in [t-C, t+C]$ we have (provided $n\ge 2$) that 
\[
\begin{split}
\|T(u, s)v\| &=\|T(u, u-C)T(u-C, s)v\| \\
&\le \theta \sup \{\|T(r, u-C)T(u-C, s)v\|: \ r\in [u-C, u+C]\} \\
&\le \theta \sup \{\|T(r, s)v\|: \ r\in [t-2C, t+2C]\}.
\end{split}
\]
Therefore,
\[
\|T(t, s)v\| \le \theta^2 \sup \{\|T(r, s)v\|: \ r\in [t-2C, t+2C]\},
\]
if $n\ge 2$.
Iterating, 
\[
\begin{split}
\|T(t, s)v\| &\le \theta^n \sup \{\|T(r, s)v\|: \ r\in [t-nC, t+nC] \} \\
&\le \theta^n \sup \{\|T(r, s)v\|: \ r\ge s\} \\
&\le D\theta^n \|v\|,
\end{split}
\]
where in the last step we used~\eqref{459new}. This easily implies that the conclusion of the lemma holds.
\end{proof}
In the sequel, we will assume that $d\in \mathbb N$. Indeed, if $d=0$ then it follows from Lemma~\ref{LEM:s} that $\mathcal T$ admits an exponential dichotomy with respect to projections $P(t)=\Id$, $t\in J$.

For $v\in Z$ with $\|v\|=1$, the map $t\mapsto \|T(t, 0)v\|$ is unbounded (as otherwise we would have $v\in\mathcal S(0)\cap Z =\{0\}$). For $n\in \mathbb N$, we define $t_n=t_n(v)>0$ by 
\[
\|T(t_n, 0)v\|=\theta^{-n}D \quad \text{and} \quad \|T(t, 0)v\|<\theta^{-n}D \ \text{for $0\le t<t_n$.}
\]
Note that $C<t_1$ (by~\eqref{HBG}). In addition, we clearly have that $t_n<t_{n+1}$ for $n\in \mathbb N$. We also claim that $t_{n+1}\le t_n+C$ for each $n\in \mathbb N$. Indeed, \eqref{531-1} gives that 
\begin{equation}\label{853}
\begin{split}
\|T(t_n, 0)v\| &=\|T(t_n, t_n-C)T(t_n-C, 0)v\|\\
&\le \theta \sup \{\|T(u, t_n-C)T(t_n-C, 0)v\|: \ u\in [t_n-C, t_n+C]\} \\
&\le \theta \sup \{\|T(u, 0)v\|: \ u \in [0, t_n+C]\} \\
&=\theta \max \{\|T(u, 0)v\|: \ u\in [0, t_n+C]\}.
\end{split}
\end{equation}
If $t_n+C<t_{n+1}$, we would have that 
\[
\max \{ \|T(u, 0)v \|: \ u \in [0, t_n+C]\}< \theta^{-(n+1)}D,
\]
and thus (by~\eqref{853}) $\|T(t_n, 0)v\| <\theta^{-n}D$, which yields a contradiction. As in~\cite[p.15]{Coppel}, this implies that 
\begin{equation}\label{interm}
\|T(t, 0)v\| \le D\theta^{-1}e^{(s-t)\frac{\log \theta}{ C}}\|T(s, 0)v\| \quad \text{for $t_1 \le t\le s<\infty$.}
\end{equation}

We now claim that $v\mapsto t_1(v)$ is bounded. Otherwise, there is a sequence $(v_n)_n\subset Z$ with $\|v_n\|=1$ such that $t_1(v_n)\to \infty$. By the compactness of the unit sphere in $Z$ (recall that $Z$ is finite-dimensional), we can assume without loss of generality that $v_n\to v$, where $v\in Z$ and $\|v\|=1$.  Then 
\[
T(t, 0)v_n\to T(t, 0)v \quad \text{for $t\in J$.}
\]
Since
\[
\|T(t, 0)v_n\| <\theta^{-1}D \quad \text{for $n\in \mathbb N$ and $0\le t<t_1(v_n)$,}
\]
and by letting $n\to \infty$ we would get that 
\[
\|T(t, 0)v\| \le \theta^{-1}D \quad \text{for $t\ge 0$,}
\]
which yields a contradiction. Set
\[
T_1:=\sup_{v\in Z, \|v\|=1}t_1(v)\in (0, \infty).
\]
Then, \eqref{interm} gives that 
\[
\|T(t, 0)v\| \le D\theta^{-1}e^{(s-t)\frac{\log \theta}{C} }\|T(s, 0)v\| \quad \text{for $T_1\le t\le s<\infty$ and $v\in Z$.}
\]
This easily implies that 
\begin{equation}\label{dich-1}
\|T(t, s)\| \le D\theta^{-1}e^{ (s-t)\frac{\log \theta}{C} } \quad \text{for $T_1\le t\le s$.}
\end{equation}
We now aim to extend~\eqref{dich-1} to the entire $J$. By Lemma~\ref{backward}, \[\sup_{t\in [0, T_1]}\|T(t, T_1)v\|<\infty \quad  \text{for all $v\in Z(T_1)$.}\]
Thus, by the Banach-Steinhaus theorem, we have that 
\begin{equation}\label{7:10}
\bar K:=\sup_{t\in [0, T_1]}\|T(t, T_1)\|<\infty. 
\end{equation}
It follows from~\eqref{dich-1} and~\eqref{7:10} that 
\begin{equation}\label{7:39}
\|T(t, s)\| \le  \bar KD\theta^{-1}e^{-\frac{T_1}{C} \log \theta}e^{(s-t)\frac{\log \theta}{C}} \quad \text{for $t\le T_1\le s<\infty$,}
\end{equation}
since $\|T(t, s)\| \le \|T(t, T_1)\| \cdot \|T(T_1, s)\|$.
Finally, \eqref{HBG} and~\eqref{7:10} give that 
\begin{equation}\label{7:37}
\|T(t, s)\| \le \|T(t, T_1)\|\cdot \|T(T_1, s)\|\le  \bar KK e^{\omega T_1} \quad \text{for $0\le t\le s\le T_1$.}
\end{equation}
It follows from~\eqref{dich-1}, \eqref{7:39} and~\eqref{7:37} that there exist $B, \alpha'>0$ such that 
\[
\|T(t, s)\| \le B'e^{-\alpha'(s-t)}, \quad \text{for $0\le t \le s$.}
\]
The above estimate together with Lemma~\ref{LEM:s} and~\cite[Lemma 4.2]{MRS} (which can be applied since $\mathcal T$ admits bounded growth) implies that $\mathcal T$ admits an exponential dichotomy.

\end{proof}

\begin{remark}
Let us compare Theorem~\ref{T1} with~\cite[Theorem 1]{D}. 
We note that in~\cite[Theorem 1]{D}, it was assumed only that $d=\codim \mathcal S(0)<\infty$. Moreover, \cite[Example 1]{D} shows that this condition cannot be omitted. For invertible evolution families this immediately implies that $\codim \mathcal S(s)=d$ for each $s\ge 0$. 
\end{remark}

\begin{example}\label{Ex1}
We now show that the conclusion of Theorem~\ref{T1} does not hold under the weaker requirement that $\codim \mathcal S(0)<\infty$.

Let $H=L^2([0, \infty))$ be the Hilbert space of all square-integrable functions $f\colon [0, \infty)\to \mathbb C$ with respect to the Lebesgue measure. We define a map $c\colon \mathbb R\to \mathbb R$ by 
\[
c(t)=\begin{cases}
1 & t<0; \\
-1 & t\ge 0.
\end{cases}
\]
For $t\ge s\ge 0$,  we define $U(t, s)\colon H\to H$ by 
\[
(U(t, s)f)(\tau)=\begin{cases}
0 & 0\le \tau <t-s; \\
e^{(t-s)c(\tau-t)}f(\tau-t+s) & \tau \ge t-s.
\end{cases}
\]
It is not difficult to verify that $\{U(t, s)\}_{t\ge s \ge 0}$ is an evolution family on $H$.

Set $X=H\oplus \mathbb R$ equipped with the norm 
\[
\|(f, y)\|=\left (\|f\|_H^2+|y|^2\right )^{1/2}, \quad (f, y)\in X.
\]
Moreover, for $t\ge s\ge 0$ we define $T(t, s)\colon X\to X$ by 
\[
T(t, s)(f, y)=(U(t, s)f, e^{t-s}y), \quad (f, y)\in X.
\]
Then $\mathcal T=\{T(t, s)\}_{t\ge s \ge 0}$ is an evolution family on $X$.

We first claim that $\mathcal T$ is exponentially expansive. To this end, for $s\ge 0$ we define $P(s)\colon X\to X$ by
\[
P(s)(f, y)=(\mathbf 1_{[s, \infty)}f, 0) \quad \text{for $(f, y)\in X$,}
\]
where $\mathbf 1_{[s, \infty)}$ denotes the characteristic function of $[s, \infty)$. Then $P(s)$ is  a projection such that $\|P(s)\|=1$ for every $s\ge 0$. A straightforward computation shows that 
\begin{equation}\label{fs1}
P(t)T(t, s)=T(t, s)P(s), \quad t\ge s\ge 0.
\end{equation}
Moreover, for $t\ge s\ge 0$ and $(f, y)\in X$ we have 
\[
\begin{split}
\|T(t, s)P(s)(f, y)\| &=\|P(t)T(t, s)P(s)(f, y)\| \\
&=\|P(t)T(t, s)(\mathbf 1_{[s, \infty)}f, 0)\| \\
&=\|(\mathbf 1_{[t, \infty)}U(t, s)(\mathbf 1_{[s, \infty)}f), 0)\| \\
&=\left (\int_t^\infty e^{2(t-s)c(\tau-t)}|f(\tau-t+s)|^2\, d\tau \right )^{1/2} \\
&=\left (\int_t^\infty e^{-2(t-s)}|f(\tau-t+s)|^2\, d\tau \right )^{1/2} \\
&\le e^{-(t-s)}\|f\|_H \\
&\le e^{-(t-s)}\|(f, y)\|.
\end{split}
\]
We conclude that
\begin{equation}\label{fs2}
\|T(t, s)P(s)\| \le e^{-(t-s)}, \quad t\ge s\ge 0.
\end{equation}
Similarly, 
\begin{equation}\label{fs3}
\|T(t, s)(f, y)\| =e^{t-s}\|(f, y)\|, \quad t\ge s\ge 0, \ (f, y)\in \Ker P(s).
\end{equation}
It follows from~\eqref{fs1}, \eqref{fs2} and~\eqref{fs3} that for $[a, b]\subset [0, \infty)$ and $t\in [a, b]$,
\[
\begin{split}
\|T(t, a)(f, y)\| &\le \|T(t, a)P(a)(f, y)\|+\|T(t, a)(\Id-P(a))(f, y)\| \\
&=\|T(t, a)P(a)(f, y)\|+\|(\Id-P(t))T(t, a)(f, y)\| \\
&\le e^{-(t-a)}\|(f, y)\|+e^{-(b-t)}\|(\Id-P(b))T(b, a)(f, y)\| \\
&\le e^{-(t-a)}\|(f, y)\|+2e^{-(b-t)}\|T(b, a)(f, y)\|,
\end{split}
\]
for every $(f, y)\in X$. This shows that $\mathcal T$ is exponentially expansive.

Note that it follows from~\eqref{fs2} and~\eqref{fs3} (applied to $s=0$) that 
\[
\mathcal S(0)=\Ima P(0)=H\oplus \{0\}.
\]
In particular, $\codim \mathcal S(0)=1<\infty$. For $s>0$, 
\[
\codim \mathcal S(s)=\codim \Ima P(s)=\dim \Ker P(s)=\dim (L^2([0, s))\oplus \mathbb R)=\infty,
\]
where we have identified the subspace of functions $f\in L^2[0, \infty)$ such that $\textbf{1}_{[s, \infty)}f=0$ with $L^2([0, s))$. This implies that $\mathcal T$ does not admit an exponential dichotomy.
\end{example}

Next, we have the version of~\cite[Theorem 2]{D} for noninvertible evolution families.
\begin{theorem}\label{T2}
Let $J=\mathbb R^+$ and $h\colon [0, \infty)\to (0, \infty)$ a growth rate. Furthermore, let  $\mathcal T=\{T(t, s): \ t\ge s, \ t, s\in J\}\subset \mathcal B(X)$ be an evolution family that exhibits $h$-bounded growth. Finally,  suppose that there exists $d\in \mathbb N_0$ such that~\eqref{codim} holds where $\mathcal S(s)$ is given by~\eqref{SS}.  Then the following statements are equivalent:
\begin{enumerate}
\item[(a)] $\mathcal T$ admits an $h$-dichotomy;
\item[(b)] $\mathcal T$ is $h$-expansive;
\item[(c)] $\mathcal T$ is uniformly $h$-noncritical.
\end{enumerate}

\end{theorem}

\begin{proof}
Let $\mathcal T_h$ be given by Theorem~\ref{TR}. It follows from~\cite[Proposition 2]{D} that $\mathcal T_h$ has bounded growth. In addition, we have the following:
\begin{itemize}
\item $\mathcal T$ has $h$-dichotomy if and only if $\mathcal T_h$ admits exponential dichotomy (see Theorem~\ref{TR});
\item $\mathcal T$ is $h$-expansive if and only if $\mathcal T_h$ is exponentially expansive (see Proposition~\ref{TR1});
\item $\mathcal T$ is uniformly $h$-noncritical if and only if $\mathcal T_h$ is uniformly noncritical (see Proposition~\ref{TR2});
\item for $s\ge 0$,
\[
\left \{v\in X: \ \sup_{t\ge s}\|T_h(t, s)v\|<\infty \right \}=\mathcal S(h^{-1}(e^s)), \quad s\in J.
\]
\end{itemize}
It remains to apply Theorem~\ref{T1}.
\end{proof}

\section{Dichotomies over $\mathbb R^-$}
We begin with an example which illustrates that Theorem~\ref{T1} does not have an immediate counterpart for the case of evolution families on $\mathbb R^-$. 

In the case of an evolution family $\{T(t, s)\}_{t\ge s, t, s\in J}$ over $J=\mathbb R^-$ on $X$, instead of subspaces $\mathcal S(s)$ introduced in~\eqref{SS}, it is natural to consider for each $s\le 0$, the subspace $\mathcal U(s)$ consisting of all $v\in X$ with the property that there is a bounded map $x\colon (-\infty, s]\to X$ such that $x(s)=v$ and $x(t)=T(t, r)x(r)$ for $r\le t\le s$.
The reason for this is that if $\mathcal T$ admits an exponential dichotomy with respect to projections $P(t)$, $t\in J$, then 
\[
\Ker P(s)=\mathcal U(s), \quad s\in J.
\]
Moreover, in this case $T(t, s)\rvert_{\mathcal U(s)}\colon \mathcal U(s)\to \mathcal U(t)$ is an isomorphism  for $t\ge s$, and consequently the natural counterpart of~\eqref{codim} would be to require that for some $d\in \mathbb N_0$,
\begin{equation}\label{dim}
\dim \mathcal U(s)=d, \quad s\in J.
\end{equation}

\begin{example}
Let $X=L^2((0, \infty))$. For $t\ge 0$, we define $S(t)\colon X\to X$ by 
\begin{equation}\label{St}
(S(t)f)(\tau)=\begin{cases}
0 & \tau <t; \\
f(\tau-t) & t\le \tau,
\end{cases}
\end{equation}
for $\tau >0$ and $f\in X$. Then $\{S(t)\}_{t\ge 0}$ is a strongly continuous semigroup of operators. We define an evolution family 
$\mathcal T=\{T(t, s)\}_{t\ge s, t, s\in J}$ over $J=\mathbb R^-$ on $X$ by 
\[
T(t, s)=e^{t-s}S(t-s), \quad 0\ge t\ge s.
\]
Then
\begin{equation}\label{3:25}
\|T(t, s)f\|=e^{t-s}\|f\|, \quad 0\ge t\ge s, \ f\in X.
\end{equation}
In particular, $\mathcal T$ has bounded growth.  For every $[a, b]\subset J$ and $a\le t\le b$ we have (using~\eqref{3:25})
\[
\|T(t, a)f\|=e^{t-a}\|f\|=e^{-(b-t)}e^{b-a}\|f\|=e^{-(b-t)}\|T(b, a)f\|,
\]
for every $f\in X$. We conclude that~\eqref{HEXP} holds with $L=\beta=1$ and $h(t)=e^t$. Thus, $\mathcal T$ is exponentially expansive.

We now claim that~\eqref{dim} holds with $d=0$. To this end, we fix $s\in J$ and take $f\in \mathcal U(s)$. In particular, for every $r\le s$ there is $g_r\in X$ such that 
\[
f=T(s, r)g_r=e^{s-r}S(s-r)g_r,
\]
which together with~\eqref{St} implies that $f$ is supported on $[s-r, \infty)$ for every $r\le s$. Hence, $f=0$ and the claim follows.

Finally, we note that $\mathcal T$ does not admit an exponential dichotomy. Indeed, otherwise,   since $\mathcal U(s)=\{0\}$ for $s\le 0$, we conclude that the associated projections $P(t)$, $t\in J$ are given by $P(t)=\Id$ for $t\in J$. Consequently, there are $D, \lambda >0$ such that 
\[
\|T(t, s)\|\le De^{-\lambda (t-s)}, \quad  0\ge t\ge s.
\]
Taking into account~\eqref{3:25}, we get that
\[
e^{t-s} \le De^{-\lambda (t-s)}, \quad 0\ge t\ge s.
\]
In particular, 
\[
e^{-s} \le De^{\lambda s}, \quad s\le 0.
\]
Letting $s\to -\infty$ yields a contradiction.
\end{example}

In the sequel, we will focus on the case when $X$ is a Hilbert space. By $A^*$ we denote the adjoint operator of $A\in \mathcal B(X)$.

The following construction will play an important role.
\begin{proposition}\label{newef}
    Let  $\mathcal T=\{T(t, s): \ t\ge s, \ t, s\in \mathbb R^-\}\subset \mathcal B(X)$ be an evolution family. For $t\ge s\ge 0$, set
    \[
    \tilde{T}(t, s):=T(-s, -t)^*.
    \]
Then $\mathcal T^*=\{\tilde{T}(t, s): \ t\ge s\ge 0\}$ is an evolution family.
\end{proposition}

\begin{proof}
Note that 
\[
\tilde{T}(t, t)=T(-t, -t)^*=\Id^*=\Id, \quad t\ge 0.
\]
For $t\ge r\ge s\ge 0$, 
\[
\begin{split}
\tilde{T}(t, r)\tilde{T}(r, s) &=T(-r, -t)^*T(-s, -r)^* \\
&=\left (T(-s, -r)T(-r, -t) \right )^* \\
&=T(-s, -t)^*\\
&=\tilde{T}(t, s).
\end{split}
\]
Finally, the strong continuity of $\tilde{\mathcal T}$ follows from~\cite[Theorem 2 (b)]{A}.
\end{proof}

\begin{proposition}\label{347prop}
        Let  $\mathcal T=\{T(t, s): \ t\ge s, \ t, s\in \mathbb R^-\}\subset \mathcal B(X)$ be an evolution family. Let $\mathcal T^*=\{\tilde T(t, s): \ t\ge s\ge 0\}$ be an evolution family given by Proposition~\ref{newef}. The following statements are equivalent:
        \begin{enumerate}
        \item[(a)] $\mathcal T$ admits an exponential dichotomy;
        \item[(b)] $\mathcal T^*$ admits an exponential dichotomy.
        \end{enumerate}
\end{proposition}

\begin{proof}
$(a)\implies (b)$ Assume that $\mathcal T$ admits an exponential dichotomy with respect to projections $P(t)$, $t\le 0$ and let $D, \lambda>0$ be such that~\eqref{HD1} and~\eqref{HD2} hold with $h(t)=e^t$. Set $\tilde{P}(t):=P(-t)^*$ for $t\ge 0$.  Using~\eqref{inv0}, we have that
\[
\begin{split}
\tilde{P}(t)\tilde{T}(t, s)=P(-t)^*T(-s, -t)^* &=(T(-s, -t)P(-t))^* \\
&=(P(-s)T(-s, -t))^* \\
&=T(-s, -t)^*P(-s)^*\\
&=\tilde{T}(t, s)\tilde P(s),
\end{split}
\]
for $t\ge s\ge 0$. Moreover, for $t\ge s\ge 0$ we have (using~\eqref{HD1}) that 
\begin{equation}\label{calc}
\begin{split}
\|\tilde{T}(t, s)\tilde P(s)\| &=\|T(-s, -t)^*P(-s)^*\| \\
&=\|(P(-s)T(-s, -t))^*\| \\
&=\|P(-s)T(-s, -t)\| \\
&=\|T(-s, -t)P(-t)\| \\
&\le De^{-\lambda (t-s)}.
\end{split}
\end{equation}
Next, we note that for $t\ge s\ge 0$
\[
\tilde{T}(t, s)\rvert_{\Ker \tilde P(s)}=T(-s, -t)^*\rvert_{\Ker P(-s)^*}\colon \Ker P(-s)^*\to \Ker P(-t)^*
\]
is  invertible by the argument in the proof of~\cite[Proposition 2.2]{RZ}. Similarly to~\eqref{calc}, \eqref{HD2} implies that 
\[
\|\tilde T(t, s)(\Id-\tilde P(s))\| \le De^{-\lambda (s-t)}, \quad 0\le t\le s.
\]
Thus $\mathcal T^*$ admits an exponential dichotomy. The converse implication can be obtained similarly.

\end{proof}

\begin{corollary}\label{cor512}
Let  $h\colon \mathbb R^-\to (0, \infty)$ be a growth rate. Let  the growth rate $k\colon \mathbb R^+\to (0, \infty)$ be given by
\begin{equation}\label{k}
k(t)=\frac{1}{h(-t)}, \quad t\ge 0.
\end{equation}
Furthermore, let $\mathcal T=\{T(t, s): \ t\ge s, \ t, s\in \mathbb R^-\}\subset \mathcal B(X)$ be an evolution family and   $\mathcal T^*=\{\tilde{T}(t, s): \ t\ge s\ge 0\}$  an evolution family given by Proposition~\ref{newef}. The following statements are equivalent:
\begin{enumerate}
\item[(a)] $\mathcal T$ admits a $h$-dichotomy;
\item[(b)] $\mathcal T^*$ admits  a $k$-dichotomy.
\end{enumerate}

\end{corollary}

\begin{proof}
By using Theorem~\ref{TR} and Proposition~\ref{347prop}, we have 
\[
\begin{split}
\text{$\mathcal T$ admits a $h$-dichotomy} &\iff \text{$\mathcal T_h$ admits an exponential dichotomy} \\
&\iff \text{$(\mathcal T_h)^*$ admits an exponential dichotomy}.
\end{split}
\]
Note that $(\mathcal T_h)^*=\{U(t, s): \ t\ge s\ge 0\}$, where
\[
\begin{split}
U(t, s)=T_h(-s, -t)^* &=T(h^{-1}(e^{-s}), h^{-1}(e^{-t}))^*\\
&=T(-k^{-1}(e^{s}), -k^{-1}(e^{t}))^*,
\end{split}
\]
for $t\ge s\ge 0$ since 
\[
k^{-1}(e^t)=-h^{-1}(e^{-t}), \quad t\ge 0.
\]
Consequently, $(\mathcal T_h)^*=(\mathcal T^*)_k$. Thus, 
\[
\text{$\mathcal T$ admits  a $h$-dichotomy} \iff \text{$(\mathcal T^*)_k$ admits an exponential dichotomy}.
\]
It remains to apply Theorem~\ref{TR} again.

\end{proof}

We now introduce a new concept of expansivity. 
\begin{definition}
    Let  $J\in \{\mathbb R, \mathbb R^+, \mathbb R^-\}$, $h\colon J\to (0, \infty)$ be a growth rate and $\mathcal T=\{T(t, s): \ t\ge s, \ t, s\in J\}\subset \mathcal B(X)$ an evolution family. We say that $\mathcal T$ is \emph{dually $h$-expansive} if there exist $L, \beta>0$ such that 
    \begin{equation}\label{dualhexp}
        \|T(b, t)^*v\| \le L\left ( \left (\frac{h(b)}{h(t)}\right )^{-\beta}\|v\|+ \left (\frac{h(t)}{h(a)}\right )^{-\beta}\|T(b, a)^*v\| \right ),
  \end{equation}
  for $v\in X$ and $t\in [a, b] \subset J$.

\end{definition}

\begin{proposition}\label{prop513}
Let  $h\colon \mathbb R^-\to (0, \infty)$ be a growth rate and $k\colon \mathbb R^+\to (0, \infty)$ given by~\eqref{k}.
 Furthermore, let $\mathcal T=\{T(t, s): \ t\ge s, \ t, s\in \mathbb R^-\}\subset \mathcal B(X)$ be an evolution family and   $\mathcal T^*=\{\tilde{T}(t, s): \ t\ge s\ge 0\}$  an evolution family given by Proposition~\ref{newef}. The following statements are equivalent:
\begin{enumerate}
\item[(a)] $\mathcal T$ is dually $h$-expansive;
\item[(b)] $\mathcal T^*$ is $k$-expansive.
\end{enumerate}
\end{proposition}

\begin{proof}
$(a)\implies (b)$ Assume that $\mathcal T$ is dually $h$-expansive and let $L, \beta>0$ be such that~\eqref{dualhexp} holds for $t\in [a, b]\subset J$. Take now $t\in [a, b]\subset \mathbb R^+$. For $v\in X$ we have 
\[
\begin{split}
\|\tilde T(t, a)v\| &=\|T(-a, -t)^*v\| \\
&\le L\left (\left (\frac{h(-a)}{h(-t)}\right )^{-\beta}\|v\|+\left (\frac{h(-t)}{h(-b)}\right )^{-\beta}\|T(-a, -b)^*v\|\right )\\
&\le L\left (\left (\frac{k(t)}{k(a)}\right )^{-\beta}\|v\|+\left (\frac{k(b)}{k(t)}\right )^{-\beta}\|\tilde{T}(b, a)v\|\right ),\\
\end{split}
\]
where we have used~\eqref{dualhexp} and $-t\in [-b, -a]\subset \mathbb R^-$. Thus, $\mathcal T^*$ is $k$-expansive. The implication $(b)\implies (a)$ follows similarly.

\end{proof}

\begin{definition}
    Let  $J\in \{\mathbb R, \mathbb R^+, \mathbb R^-\}$, $h\colon J\to (0, \infty)$ be a growth rate and $\mathcal T=\{T(t, s): \ t\ge s, \ t, s\in J\}\subset \mathcal B(X)$ an evolution family. We say that $\mathcal T$ is \emph{dually uniformly $h$-noncritical} if there exist $\theta \in (0, 1)$ and $C>0$ such that 
    \begin{equation}\label{newnon}
    \begin{split}
  & \|T(h^{-1}(e^C h(t)),  t)^*v\| \\
  &\le \theta \sup \{\|T(h^{-1}(e^C h(t)), u)^*v\|: \ h^{-1}(e^{-C}h(t))\le u \le h^{-1}(e^C h(t))\},
 \end{split}
    \end{equation}
for all $v\in X$ and $t\in J$ such that $[h^{-1}(e^{-C} h(t)), h^{-1}(e^C h(t))]\subset J$.

\end{definition}

\begin{proposition}\label{prop515}
    Let  $h\colon \mathbb R^-\to (0, \infty)$ be a growth rate. We define a growth rate $k\colon \mathbb R^+\to (0, \infty)$ by~\eqref{k}.  Furthermore, let  $\mathcal T^*=\{\tilde{T}(t, s): \ t\ge s\ge 0\}$ be an evolution family given by Proposition~\ref{newef}. The following statements are equivalent:
\begin{enumerate}
\item[(a)] $\mathcal T$ is dually uniformly $h$-noncritical;
\item[(b)] $\mathcal T^*$ is uniformly $k$-noncritical.
\end{enumerate}

\end{proposition}

\begin{proof}
$(a)\implies (b)$ Suppose that  $\mathcal T$ is dually uniformly $h$-noncritical and let $\theta \in (0, 1)$ and $C>0$ be such that~\eqref{newnon} holds for all $v\in X$ and $t\in \mathbb R^-$ such that $[h^{-1}(e^{-C} h(t)), h^{-1}(e^C h(t))]\subset \mathbb R^-$.

Take now $v\in X$ and $t\in \mathbb R^+$ such that  $[t, k^{-1}(e^{2C} k(t))]\subset \mathbb R^+$. 
It follows from~\eqref{k} that
\[
k^{-1}(e^C k(t))=-h^{-1}(e^{-C}h(-t)), \quad t\ge 0,
\]
and
\[
k^{-1}(e^{2C}k(t))=-h^{-1}(e^{-2C}h(-t)), \quad t\ge 0.
\]
Thus,
\[
\begin{split}
\|\tilde T(k^{-1}(e^C k(t)), t)v\| &=\|T(-t, -k^{-1}(e^C k(t)))^*v\| \\
&=\|T(-t, h^{-1}(e^{-C}h(-t)))^*v\| \\
&\le \theta \sup \{\|T(-t, u)^*v\|: \ h^{-1}(e^{-2C}h(-t)) \le u \le -t\} \\
&= \theta \sup \{\|T(-t, -u)^*v\|: \ t \le u \le -h^{-1}(e^{-2C}h(-t))\} \\
&= \theta \sup \{\|\tilde T(u, t)v\|: \ t \le u \le -h^{-1}(e^{-2C}h(-t))\} \\
&= \theta \sup \{\|\tilde T(u, t)v\|: \ t \le u \le k^{-1}(e^{2C}k(t))\},
\end{split}
\]
where we applied~\eqref{newnon} for $h^{-1}(e^{-C}h(-t))$ instead of $t$. Similarly one can establish the converse. 
\end{proof}

The following is a central result of this section.
\begin{theorem}
Let $h\colon (-\infty, 0] \to (0, \infty)$ be  a growth rate and       $\mathcal T=\{T(t, s): \ t\ge s, \ t, s\in \mathbb R^-\}\subset \mathcal B(X)$ an evolution family with $h$-bounded growth. Furthermore, suppose that there is $d\in \mathbb N_0$ such that 
\begin{equation}
\codim \mathcal S^*(s)=d \quad \text{for $s\le 0$,}
\end{equation}
where
\[
\mathcal S^*(s)=\left \{v\in X: \ \sup_{t\le s}\|T(s, t)^*v\|<\infty \right \}.
\]
Then the following statements are equivalent:
\begin{enumerate}
\item[(a)] $\mathcal T$  has $h$-dichotomy;
\item[(b)] $\mathcal T$ is dually $h$-expansive;
\item[(c)] $\mathcal T$ is dually uniformly $h$-noncritical.
\end{enumerate}
\end{theorem}

\begin{proof}
Let $\mathcal T^*=\{\tilde T(t, s): \ t\ge s\ge 0\}$ be an evolution family given by Proposition~\ref{newef}. Observe that 
\[
\mathcal S^*(s)=\left \{v\in X:  \ \sup_{t\ge -s}\|\tilde{T}(t, -s)v\|<\infty \right \}, \quad s\le 0.
\]
Furthermore, we note that $\mathcal T^*$ has $k$-bounded growth, where $k$ is given by~\eqref{k}. Indeed, since $\mathcal T$ has $h$-bounded growth, there are $D, \omega>0$ such that 
\[
\|T(t, s)\| \le D\left (\frac{h(t)}{h(s)}\right )^\omega, \quad 0\ge t\ge s.
\]
Consequently, for $t\ge s\ge 0$ we have 
\[
\|\tilde T(t, s)\|=\|T(-s, -t)^*\|=\|T(-s, -t)\| \le D\left (\frac{h(-s)}{h(-t)}\right )^\omega =D\left (\frac{k(t)}{k(s)}\right )^\omega.
\]
In addition, we have the following:
\begin{itemize}
\item $\mathcal T$ admits a $h$-dichotomy if and only if $\mathcal T^*$ admits a $k$-dichotomy (see Corollary~\ref{cor512});
\item $\mathcal T$ is dually $h$-expansive if and only if $\mathcal T^*$ is $k$-expansive  (see Proposition~\ref{prop513});
\item $\mathcal T$ is dually uniformly $h$-noncritical if and only if $\mathcal T^*$ is uniformly $k$-noncritical (see Proposition~\ref{prop515}).
\end{itemize}
The conclusion of the theorem now follows readily from Theorem~\ref{T2}.
\end{proof}

\vspace{7mm}

\noindent {\bf  Funding}  
This paper has been funded by European Union – NextGenerationEU-Statistical properties of random dynamical systems and other contributions to mathematical analysis and probability theory.

\vspace{3mm}

\noindent{\bf Data availability} This  manuscript has no associated data.

\section*{Declarations}

\noindent{\bf Conflict of interest}  The author declares that he has no known competing financial interests or personal
relationships that could have appeared to influence the work in this paper.

\end{document}